\documentclass[10pt]{article}
\usepackage{latexsym}
\usepackage{amssymb}
\usepackage[bookmarks=false]{hyperref}
\usepackage{ifpdf}
\begin{document}
\title{{\bf Hom-Lie-Yamaguti Structures on Hom-Leibniz Algebras}}
\author{{\bf Donatien Gaparayi} \footnote{Permanent address: Ecole Normale Sup\'erieure (E.N.S), BP 6983 Bujumbura, Burundi.}\\Institut de Math\'ematiques et de Sciences Physiques,\\ 01 BP 613-Oganla, Porto-Novo, B\'enin\\gapadona@yahoo.fr\vspace{0.30cm} \\
{\bf A. Nourou Issa} \\ D\'epartement de Math\'ematiques, Universit\'e d'Abomey-Calavi\\01 BP 4521, Cotonou, B\'enin\\ woraniss@yahoo.fr}
\date{}
\maketitle
\begin{abstract}
Every multiplicative left Hom-Leibniz algebra has a natural Hom-Lie-Yamaguti structure.\\
\\
AMS {\it Subject Classification (2010): 17A30, 17A32, 17D99}\vspace{0.30cm}
\end{abstract}
\par
{\bf 1. Introduction and statement of result} \\
\par
A {\it (left) Leibniz algebra} is an algebra $(L,\cdot)$ satisfying the identity \\
\par
$ x \cdot (y \cdot z) = (x \cdot y) \cdot z + y \cdot (x \cdot z)$. \\
\\
Leibniz algebras were introduced by J.-L. Loday \cite{Loday} (and so they are sometimes called Loday algebras) as a noncommutative analogue of Lie algebras, in the study of some topics in homological algebra  and noncommutative geometry (see also \cite{Casas}, \cite{Pirashvili}). While earlier papers on Leibniz algebras are concerned with some homological problems (see, e.g., \cite{Loday}, \cite{Pirashvili}), some structure theory of Leibniz algebras are proposed in, e.g., \cite{Omirov} and \cite{Ayupov} (see also references therein). Classification of low-dimensional Leibniz algebras could be found in, e.g., \cite{Omirov}, \cite{Felipe}, \cite{Loday}, \cite{Rakh}.
\par
One of the problems in the general theory of a given class of (binary or binary-ternary) nonassociative
algebras is the study of relationships between that class of algebras and the one of Lie algebras. In the same rule, the search of relationships between a class of nonassociative algebras and the one of Leibniz algebras is of interest (at least for constructing concrete examples of the given class of nonassociative algebras). In this setting, the existence of a Lie-Yamaguti structure on any (left) Leibniz algebra pointed out in \cite{KIN1} is a good illustration.
\par
A {\it Lie-Yamaguti algebra} is a triple $(L,[,],\{,,\})$ in which $L$ is a vector space,  $[,] : L \times L \rightarrow L$ a bilinear map and  
$\{,,\} : L \times L\times L \rightarrow L$  a trilinear map such that
\par
(LY1) $[x,y] = - [y,x],$
\par
(LY2) $\{x,y,z\} = -\{y,x,z\},$
\par
(LY3) $\circlearrowleft_{x,y,z}([[x,y],z] + \{x,y,z\}) = 0,$
\par
(LY4) $\circlearrowleft_{x,y,z} \{[x,y],z,u\} = 0,$
\par
(LY5) $\{x,y,[u,v]\} = [\{x,y,u\},v] + [u,\{x,y,v\}],$
\par
(LY6) $\{x,y,\{u,v,w\}\} = \{\{x,y,u\},v,w\} + \{u,\{x,y,v\},w\}$
\par
\hspace{4truecm} $+ \{u,v,\{x,y,w\}\},$ \\
for all $u,v,w,x,y,z,$ in $L$, where $\circlearrowleft_{x,y,z}$ denotes the sum over cyclic permutation of $x,y,z.$ Lie -Yamaguti algebras, first called ''generalized Lie triple systems'', were introduced by K. Yamaguti \cite{LY1} while giving an algebraic interpretation of the characteristic properties 
of the torsion and curvature of a homogeneous space with canonical connection (the Nomizu's connection) \cite{Nom1}. Later, M. Kikkawa \cite{KIK1} called them ''Lie triple algebras'' and the terminology of ''Lie-Yamaguti algebras'' is introduced in \cite{KIN1} to designate these algebras. As in \cite{Benito}, we write ''LY-algebras'' for Lie-Yamaguti algebras.
\par
In a left Leibniz algebra $(L,\cdot)$ if define $[x,y] := x \cdot y - y \cdot x$ (skew- symmetrization) and $\{x,y,z\} := -(x \cdot y) \cdot z$, then $(L,[,],\{,,\})$ is a LY-algebra \cite{KIN1}. In this note, we will be interested in the counterpart of this construction in the Hom-algebra setting.
\par
Roughly speaking, a Hom-type generalization of a kind of algebras is defined by twisting its defining identities by a linear self-map (the twisting map) in such a way that when the twisting map is the identity map, one recovers the original kind of algebras. The theory of Hom-algebras originated from the introduction of so-called ''Hom-Lie algebras'' in \cite{HAR1} as an abstraction of an approach to deformations of the Witt algebra and the Virasoro algebra based on $\sigma$-derivations, including $q$-derivations of the Witt algebra and the Virasoro algebra associated to $q$-difference operators. The outcome of the algebraic
view of these considerations is the introduction of ''Hom-associative algebras'' in \cite{MAK3} as a Hom-analogue of associative algebras. Hom-associative algebras are for Hom-Lie algebras what are associative algebras for Lie algebras: the commutator-algebra of a Hom-associative algebra is a Hom-Lie algebra \cite{MAK3}. As a ''noncommutative'' generalization of Hom-Lie algebras, Hom-Leibniz algebras are also defined in \cite{MAK3}, where other results regarding Hom-algebras are found (see also \cite{YAU1}).
A general method of twisting ordinary algebras into their Hom-type analogues is given in \cite{YAU2}. The reader is refered to, e.g., \cite{GAP}, \cite{NOU1},\cite{MAK1}, \cite{MAK2}, \cite{YAU3} for discussions about various Hom-type algebras.
\par
Following this line, Hom-Lie-Yamaguti algebras (Hom-LY algebras) are introduced in \cite{GAP}.  A {\it Hom-Lie-Yamaguti algebra} (Hom-LY algebra for 
short) is a quadruple $(L,[,],\{,,\},\alpha )$ in which $L$ is a vector space, $"[,]"$ a binary operation and $"\{,,\}"$ a ternary operation on $L$, and $\alpha : L \rightarrow L$ a linear map such that
\par
(HLY1) $\alpha ([x,y]) = [\alpha(x),\alpha(y)],$
\par
(HLY2) $\alpha (\{x,y,z\}) = \{\alpha(x),\alpha(y),\alpha(z)\},$ 
\par
(HLY3) $[x,y] = - [y,x],$
\par
(HLY4) $\{x,y,z\} = - \{y,x,z\},$ 
\par
(HLY5) $\circlearrowleft_{x,y,z} ([[x,y],\alpha(z)] + \{x,y,z\}) = 0,$ 
\par
(HLY6) $\circlearrowleft_{x,y,z}\{[x,y], \alpha(z),\alpha(u)\} = 0,$
\par
(HLY7) $\{\alpha(x),\alpha(y),[u,v]\} = [\{x,y,u\},\alpha^2(v)] + [\alpha^2(u),\{x,y,v\}],$ 
\par
(HLY8) $\{\alpha^2(x),\alpha^2(y),\{u,v,w\}\} = \{\{x,y,u\},\alpha^2(v),\alpha^2(w)\}$
\par
\hspace{3truecm} $+ \{\alpha^2(u),\{x,y,v\},\alpha^2(w)\} + \{\alpha^2(u),\alpha^2(v),\{x,y,w\}\},$ \\
for all $u,v,w,x,y,z \in Y.$  Thus, as mentioned above, we shall prove the following\\
\par
{\bf Theorem.} {\it Every multiplicative left Hom-Leibniz algebra has a natural Hom-Lie-Yamaguti structure.} \\
\par
The useful definitions and some facts as the characterization of the Hom-Akivis algebra associated to a given Hom-Leibniz algebra are reminded in section 2. In section 3, we prove the theorem and discuss examples of Hom-LY algebras that are constructed using the theorem above (thusly, we also construct examples of left  Hom-Leibniz algebras).
\par
All vector spaces and algebras are considered over a fixed ground field of characteristic 0.\\
\par
{\bf 2. Definitions and basic facts} \\
\par
We recall some basic notions, introduced in \cite{HAR1}, \cite{NOU1}, \cite{MAK3}, \cite{YAU1}, \cite{YAU2}, related to Hom-algebras. We also recall from
\cite{NOU2} a characterization of the Hom-Akivis algebra associated with a given Hom-Leibniz algebra.\\
\par
{\bf Definition 2.1.} (\cite{MAK3}, \cite{YAU1}) A {\it Hom-algebra} is a triple $(L, \cdot ,\alpha)$ in which $L$ is a vector space, ''$\cdot$'' a binary operation on $L$ and  $\alpha : L \rightarrow L$ is a linear map (the twisting map).
\par
A Hom-algebra $(L, \cdot,\alpha)$ is said to be {\it multiplicative} if $\alpha(x \cdot y) = \alpha(x) \cdot \alpha(y)$ (multiplicativity), for all $x,y$ in $L$.
\par
Since our result depends on multiplicativity, we assume here that all Hom-algebras are multiplicative. \\
\par
{\bf Definition 2.2.} Let $(L, \cdot ,\alpha)$ be a Hom-algebra.\\
(i) The {\it Hom-associator} \cite{MAK3} of $L$ is the trilinear map $as_\alpha : L\times L\times L \rightarrow L$ defined by \\
\par
$as_\alpha(x,y,z) = (x.y).\alpha(z) - \alpha(x).(y.z),$ \hfill (2.1)\\
\\
for all $x,y,z \in L.$ If $as_\alpha(x,y,z) = 0$ (Hom-associativity), $\forall x,y,z \in L,$ then $(L,\cdot,\alpha)$ is said to be {\it Hom-associative} \cite{MAK3}.\\
(ii) The {\it Hom-Jacobian} \cite{MAK3} of $L$ is the trilinear map $J_\alpha : L\times L\times L \rightarrow L$ defined by\\
\par
$J_\alpha (x,y,z) := \circlearrowleft_{x,y,z}(x.y).\alpha(z)$ \hfill $(2.2)$\\
\\
for all $x,y,z$ in $L$. The Hom-algebra $(L,\cdot ,\alpha)$ is called a {\it Hom-Lie algebra} \cite{HAR1} if the operation ''$\cdot$'' is anticommutative and the {\it Hom-Jacobi identity} $J_\alpha (x,y,z) = 0$ is satisfied in $(L,\cdot,\alpha)$.\\
\par
{\bf Remark 2.3.} If $\alpha = Id$ (the identity map) then (2.1) (resp. (2.2)) is just the associator (resp. the Jacobian) in $(L, \cdot ,\alpha).$ Therefore an associative (resp. a Lie) algebra could be seen as a Hom-associative (resp. Hom-Lie) algebra with the identity map as the twisting map. Also note that a not necessarily Hom-associative algebra is called a {\it non-Hom-associative} algebra in \cite{NOU1} in analogy with the case of not necessarily associative algebras (the terminologies of ''Hom-nonassociative algebras'' or ''nonassociative Hom-algebras'' are also used in \cite{MAK2}, \cite{YAU1} for that type of Hom-algebras).
\par
As for Lie algebras, Hom-Lie algebras have a ''noncommutative'' generalization as Hom-Leibniz algebras. \\
\par
{\bf Definition 2.4.} (\cite{MAK3}) A {\it (left) Hom-Leibniz algebra} is a Hom-algebra $(L,\cdot ,\alpha)$ satisfying the {\it (left) Hom-Leibniz identity}\\
\par
$\alpha(x) \cdot (y \cdot z) = (x \cdot y) \cdot \alpha(z) + \alpha(y) \cdot (x \cdot z) $ \hfill (2.3)\\
\par
for all $x,y,z$ in $L$.\\
\par
{\bf Remark 2.5.} If $\alpha = Id$ in Definition 2.4, then $(L, \cdot ,\alpha)$ reduces to a (left) Leibniz algebra $(L, \cdot).$ Moreover, as for Leibniz algebras \cite{Loday}, if the operation of a given Hom-Leibniz algebra  $(L,\cdot,\alpha)$ is skew-symmetric (i.e. anticommutative), then $(L,\cdot,\alpha)$ turns out to 
be a Hom-Lie algebra (see \cite{MAK3}). We also observe that the original definition of a Hom-Leibniz algebra \cite{MAK3} is related to the identity (the ''right'' Hom-Leibniz identity)\\
\par
$(x\cdot y) \cdot \alpha(z) = (x\cdot z)\cdot \alpha(y) + \alpha(x)\cdot (y\cdot z).$\\
\\
Moreover, given a linear self-map $\alpha$ of $L$, every Leibniz algebra $(L,\cdot)$ can be twisted into a Hom-Leibniz algebra $(L,\dot{\alpha},\alpha)$ with ``$\dot{\alpha}$`` defined by $x \; \dot{\alpha} \; y = \alpha (x \cdot y)$ for all $x,y$ in $L$ (\cite{YAU2}). \\
\par
The extension to binary-ternary algebras of twisting identities of algebras is considered in \cite{NOU1}. This led to the introduction of the class of ''Hom-Akivis algebras'' as a twisted version of Akivis algebras introduced by M.A. Akivis (see \cite{Akivis1} and references therein).
\par
Akivis algebras (first called ''$W$-algebras '' \cite{Akivis1}) arose in the differential geometry of differentiable webs, and also as tangent algebras to local differentiable quasigroups. The terminology of ''Akivis algebras'' is introduced in \cite{Hofmann1}.\\
\par
{\bf Definition 2.6.} (\cite{NOU1}) A {\it Hom-Akivis algebra } is a quadruple $(L,[,],[,,],\alpha)$ such that $L$ is a vector space, ''$[,]$'' is a 
skew-symmetric binary operation on $L$, ''$[,,]$'' a ternary operation on $L,$  $\alpha : L \rightarrow L$ a linear map, and that the {\it Hom-Akivis identity}\\
\par
$J_\alpha (x,y,z) = \circlearrowleft_{x,y,z}[x,y,z] - \circlearrowleft_{x,y,z}[y,x,z]$ \hfill (2.4)\\
\\
holds for all $x,y,z$ in $L.$\\
\par
Observe that for $\alpha = Id$ the Hom-Akivis identity $(2.4)$ reduces to the {\it Akivis identity}\\
\par
$J (x,y,z) = \circlearrowleft_{x,y,z}[x,y,z] - \circlearrowleft_{x,y,z}[y,x,z]$\\
\\
which defines {\it Akivis algebras.} It is shown \cite{NOU1} that every Akivis algebra with a linear self-map is twisted into a Hom-Akivis algebra and that every non-Hom-associative algebra with a linear self-map is a Hom-Akivis algebra with respect to the skew-symmetrization $[x,y] = xy - yx$ and Hom-associator
$[x,y,z] = as_{\alpha}(x,y,z).$
\par
In terms of Hom-associators, the identity (2.3) has the form \\
\par
$as_{\alpha}(x,y,z) = -\alpha(y).(x.z)$. \hfill (2.5)\\
\par
Thus the operations of the Hom-Akivis algebra associated with the Hom-Leibniz algebra $(L,\cdot,\alpha)$ are the skew-symmetrization and (2.5). Then the Hom-Akivis identity (2.4) takes the form\\
\par
$\circlearrowleft_{x,y,z}[[x,y],\alpha(z)] = \circlearrowleft_{x,y,z}as_{\alpha}(x,y,z) - \circlearrowleft_{x,y,z}as_{\alpha}(y,x,z)$\\
\\
that is, by (2.5) and (2.3),\\
\par
$\circlearrowleft_{x,y,z}[[x,y],\alpha(z)]=\circlearrowleft_{x,y,z}(x \cdot y)\cdot \alpha(z)$.\hfill (2.6)\\
\\
The considerations above will be used in the next section in the proof of the theorem.\\
\par
{\bf 3. Proof of the Theorem. Examples} \\
\par
In this section, we settle down in the proof of our claim, i.e. the existence of a Hom-LY structure on any (multiplicative) left Hom-Leibniz algebra. This proof is based on a specific ternary operation that can be considered on a given Hom-Leibniz algebra (this product is the Hom-analogue of the ternary
operation considered in \cite{KIN1} on a left Leibniz algebra $L$ that produces, along with the skew-symmetrization, a LY structure on $L$). Also note that our proof below essentially relies on some properties characterizing Hom-Leibniz algebras, obtained in  \cite{NOU2}. We conclude, as an illustration of our result, by some constructions of Hom-LY algebras from twisted Leibniz algebras (incidentally, this produces examples of left Hom-Leibniz algebras).
\par
In a left Hom-Leibniz algebra $(L,\cdot ,\alpha)$ consider the skew-symmetrization\\
\par
$[x,y] := x \cdot y - y \cdot x$ \\
\\
for all $x,y$ in $L.$ Then, from \cite{NOU2}, we know that \\
\par
$(x \cdot y + y \cdot x) \cdot \alpha(z) = 0$, \hfill (3.1)
\par
$\alpha(x) \cdot [y,z] = [(x \cdot y),\alpha(z)] + [\alpha(y),(x \cdot z)].$ \hfill (3.2)\\
\par
If consider the left translations $\Lambda_{a} b := a \cdot b$ in $(L, \cdot ,\alpha),$ then the identities (2.3) and (3.2) can be written respectively as\\
\par
$\Lambda_{\alpha(x)}(y \cdot z) = (\Lambda_{x} y)\cdot \alpha(z) + \alpha(y) \cdot (\Lambda_{x} z),$ \hfill (3.3)
\par
$\Lambda_{\alpha(x)}[y,z] = [\Lambda_{x}y,\alpha(z)] + [\alpha(y),\Lambda_{x}z].$ \hfill (3.4)\\
\\
{\it Proof of the Theorem.}\\
\par
In $(L, \cdot ,\alpha)$ consider the following ternary operation: \\
\par
$\{x,y,z\} := as_\alpha(y,x,z) - as_\alpha(x,y,z)$ \hfill (3.5)\\
\\
for all $x,y,z$ in $L$. Then (3.5), (2.3) and (2.5) imply \\
\par
$\{x,y,z\} = - (x \cdot y) \cdot \alpha(z)$. \hfill (3.6)\\
\\
Moreover, we have
\begin{eqnarray*}
 [x,y] \cdot \alpha(z) & = & (x \cdot y - y \cdot x) \cdot \alpha(z)\\
& = & 2 (x \cdot y) \cdot \alpha(z)\quad \mbox{(by (3.1))}\\ 
& = & -2 \{x,y,z\} \quad \mbox{(see (3.6))}
\end{eqnarray*}
so that\\
\par
$\{x,y,z\} = -\frac{1}{2}[x,y] \cdot \alpha(z)$. \hfill (3.7)\\
\\
Thus (3.5), (3.6) and (3.7) are different expressions of the operation ''$\{,,\}$'' that are for use in what follows.
\par
Now we proceed to verify the validity on $(L, \cdot ,\alpha)$ of the set of identities (HLY1)-(HLY8).
\par
The multiplicativity of $(L, \cdot ,\alpha)$ implies (HLY1) and (HLY2) while (HLY3) is the skew-symmetrization and (HLY4) clearly follows from (3.5) (or (3.7)). Next, observe that (HLY5) is just the Hom-Akivis identity (2.6) for $(L, \cdot ,\alpha).$
\par
Consider now $\circlearrowleft_{(x,y,z)}\{[x,y],\alpha(z),\alpha(u)\}$. 
Then
\begin{eqnarray*}
\circlearrowleft_{x,y,z}\{[x,y], \alpha(z),\alpha(u)\}
& = & \circlearrowleft_{x,y,z}-([x,y] \cdot \alpha(z)) \cdot \alpha^2(u) \quad(\mbox{by (3.6)})\\
& = & 2 (\circlearrowleft_{x,y,z}\{x,y,z\}) \cdot \alpha^2(u) \quad(\mbox{by (3.7)})\\
& = & -2 ((x \cdot y) \cdot \alpha(z) + (y \cdot z) \cdot \alpha(x) + (z \cdot x) \cdot \alpha(y)) \cdot \alpha^2(u)\\
& = & -2 (\alpha(x) \cdot (y \cdot z) - \alpha(y) \cdot (x \cdot z) + (y \cdot z) \cdot \alpha(x)\\
& & +(z \cdot x) \cdot \alpha(y)) \cdot \alpha^2(u)\quad(\mbox{by (2.3)})\\
& = & -2 (\alpha(x) \cdot (y \cdot z) + (y \cdot z) \cdot \alpha(x)) \cdot \alpha^2(u)\\
&& -2 (- \alpha(y) \cdot (x \cdot z) + (z \cdot x) \cdot \alpha(y)) \cdot \alpha^2(u)\\
& = & -2 (- \alpha(y) \cdot (x \cdot z) + (z \cdot x) \cdot \alpha(y)) \cdot \alpha^2(u) \quad(\mbox{by (3.1)})\\
& = & -2 (- \alpha(y) \cdot (x \cdot z) - (x \cdot z) \cdot \alpha(y)) \cdot \alpha^2(u) \quad(\mbox{by (3.1)})\\
& = & 2 (\alpha(y) \cdot (x \cdot z) + (x \cdot z) \cdot \alpha(y)) \cdot \alpha^2(u)\\
& = & 0 \quad(\mbox{by (3.1)})
\end{eqnarray*}
so that we get (HLY6). Next
\begin{eqnarray*}
\{\alpha(x),\alpha(y),[u,v]\} & = & - \alpha(x \cdot y) \cdot \alpha([u,v])\quad(\mbox{by (3.6) and multiplicativity})\\ 
& = & \Lambda_{-\alpha(x \cdot y)}[\alpha(u),\alpha(v)] \\
& = &[\Lambda_{-x \cdot y}\alpha(u),\alpha^2(v)] + [\alpha^2(u),\Lambda_{-x \cdot y}\alpha(v)] \quad(\mbox{by (3.4)})\\
& = & [\{x,y,u\},\alpha^2(v)]  + [\alpha^2(u),\{x,y,v\}]\quad(\mbox{by (3.6)})
\end{eqnarray*}
which is (HLY7). Finally, we compute
\begin{eqnarray*}
&&\{\{x,y,u\},\alpha^2(v),\alpha^2(w)\} + \{\alpha^2(u),\{x,y,v\},\alpha^2(w)\}\\&& +\{\alpha^2(u),\alpha^2(v),\{x,y,w\}\}\\
& = & \{-\Lambda_{x \cdot y}\alpha(u),\alpha^2(v),\alpha^2(w)\} + \{\alpha^2(u),-\Lambda_{x \cdot y} \alpha(v),\alpha^2(w)\} \\
&& + \{\alpha^2(u),\alpha^2(v),-\Lambda_{x \cdot y}\alpha(w)\}\\
& = & -((-\Lambda_{x \cdot y}\alpha(u)) \cdot \alpha^2(v)) \cdot \alpha^3(w) - (\alpha^2(u) \cdot (-\Lambda_{x \cdot y}\alpha(v))) \cdot \alpha^3(w) \\
&&- (\alpha^2(u) \cdot \alpha^2(v))\alpha(-\Lambda_{x \cdot y}\alpha(w))\\
& = & (\Lambda_{\alpha(x \cdot y)}\alpha(u \cdot v)) \cdot \alpha^3(w) + \alpha^2(u \cdot v) \cdot \Lambda_{\alpha(x \cdot y)}\alpha^2(w)\quad(\mbox{by (3.3) and multiplicativity})\\
& = & \Lambda_{\alpha^2(x \cdot y)}(\alpha(u \cdot v) \cdot \alpha^2(w)) \quad(\mbox{by (3.3)})\\
& = & -\alpha^2(x \cdot y) \cdot (-\alpha(u \cdot v) \cdot \alpha^2(w))\\
& = & -(\alpha^2(x) \cdot \alpha^2(y)) \cdot \alpha(-(u \cdot v) \cdot \alpha(w))\quad(\mbox{by multiplicativity})\\
& = & \{\alpha^2(x),\alpha^2(y),\{u,v,w\}\} \quad(\mbox{by (3.6)}).
\end{eqnarray*}
Therefore $(L,[,],\{,,\},\alpha)$ is a Hom-LY algebra.
This completes the proof.
\hfill $\square$ \\
\par
{\bf Remark 3.1.} If set $\alpha = Id$ in (3.6), then we recover the ternary operation defined in \cite{KIN1} in the proof of the existence of a natural LY structure on any left Leibniz algebra (see section 1). Therefore, although with a quite different scheme of proof, our result here is an 
$\alpha$-twisted version of the one in \cite{KIN1}. The untwisted version of the proof proposed here could be found in \cite{NOU3}, where the result of \cite{KIN1} is considered again but via Akivis algebras.\\
\par
We now discuss examples of Hom-LY algebras that can be constructed using the theorem. Examples of Hom-LY algebras constructed from LY algebras can be found in \cite{GAP}.\\
\par
{\bf Example 3.2.} Let $(L, \cdot , \alpha)$ be an anticommutative multiplicative left Hom-Leibniz algebra. Then $(L, \cdot , \alpha)$ is a multiplicative Hom-Lie algebra (\cite{MAK3}). If define on $L$ a ternary operation  ``$\{, , \}$'' by (3.6), then one checks that $(L, \cdot , \{ , , \}, \alpha)$ is a Hom-LY algebra. \\
\par
In the following examples, the unspecified products are regarded as zero. \\
\par
{\bf Example 3.3.} Let $(L, \cdot )$ be a 3-dimensional complex algebra defined by \\
\par
$e_2 \cdot e_3 = e_2$, $e_3 \cdot e_1 = {\lambda}e_1$, $e_3 \cdot e_2 = -e_2$, $e_3 \cdot e_3 = e_1$ $(\lambda \in \mathbb{C})$.\\
\\
Then $(L, \cdot )$ is a (solvable) complex left Leibniz algebra (\cite{Rakh}, Table 3, algebra $L_4$). Now a few computation shows that all the linear self-maps $\alpha$ of $L$ given by\\
\par
$\alpha (e_1) = (a \lambda + 1)e_1$, $\alpha (e_2) = b e_2$, $\alpha (e_3) = a e_1 + e_2 + e_3$\\
\\
are endomorphisms of $(L, \cdot )$, where $a,b, \lambda \in \mathbb{C}$. By a result in \cite{YAU2}, we define on $L$ an operation ``$\dot{\alpha}$'' by \\
\par
$e_2 \; \dot{\alpha} \; e_3 := \alpha (e_2 \cdot e_3) = b e_2$,
\par
$e_3 \; \dot{\alpha} \; e_1 := \alpha (e_3 \cdot e_1) = \lambda (a \lambda + 1)e_1$,
\par
$e_3 \; \dot{\alpha} \; e_2 := \alpha (e_3 \cdot e_2) = - b e_2$,
\par
$e_3 \; \dot{\alpha} \; e_3 := \alpha (e_3 \cdot e_3) = (a \lambda + 1)e_1$ \\
\\
to get a multiplicative left Hom-Leibniz algebra $L_{\alpha} := (L, \dot{\alpha}, \alpha)$. Therefore, according to the theorem, the Hom-LY algebra $(L, [,] , \{ , , \}, \alpha)$ corresponding to $L_{\alpha}$ is defined by \\
\par
$[e_1 , e_3] = -\lambda (a \lambda + 1)e_1 \;$ ($=-[e_3 , e_1]$),
\par
$[e_2 , e_3] = 2b e_2 \;$ ($=-[e_3, e_2]$),
\par
$\{ e_3, e_2, e_3 \} = b^2 e_2 \;$ ($=-\{ e_2, e_3, e_3 \}$) \\
\\
with $a,b, \lambda \in \mathbb{C}$.\\
\par
{\bf Example 3.4.} Let $(L, \cdot )$ be a 4-dimensional complex algebra defined by\\
\par
$e_1 \cdot e_1 = e_4$, $e_1 \cdot e_2 = e_3$, $e_1 \cdot e_3 = e_4$,
\par
$e_2 \cdot e_1 = -e_3$,
\par
$e_3 \cdot e_1 = -e_4$.\\
\\
Then $(L, \cdot )$ is a (nilpotent) complex left Leibniz algebra (\cite{Omirov}, Theorem 3.2, algebra $\mathfrak{R}_7$). If define a linear self-map of $L$ by \\
\par
$\alpha (e_1) = e_1 + e_2 + e_3 + e_4$,
\par
$\alpha (e_2) = e_2 + e_3 + e_4$,
\par
$\alpha (e_3) = e_3 + e_4$,
\par
$\alpha (e_4) = e_4$,\\
\\
then $\alpha$ is verified to be an endomorphism of $(L, \cdot )$. Next, as in the example above, we may define on $L$ an operation ``$\dot{\alpha}$'' by \\
\par
$e_1 \; \dot{\alpha} \; e_1 = e_4$, $e_1 \; \dot{\alpha} \; e_2 = e_3 + e_4$, $e_1 \; \dot{\alpha} \; e_3 = e_4$,
\par
$e_2 \; \dot{\alpha} \; e_1 = -e_3 - e_4$,
\par
$e_3 \; \dot{\alpha} \; e_1 = -e_4$\\
\\
and then $L_{\alpha} := (L, \dot{\alpha}, \alpha)$ is a multiplicative left Hom-Leibniz algebra. Now, applying the theorem, we see that the Hom-LY algebra $(L, [,] , \{ , , \}, \alpha)$ corresponding to $L_{\alpha}$ is constructed by defining\\
\par
$[e_1 , e_2] = 2(e_3 + e_4) \;$ ($=-[e_2, e_1]$), $\;[e_1 , e_3] = 2e_4 \;$ ($=-[e_3, e_1]$),
\par
$\{ e_1, e_2, e_1 \} = e_4 \;$ ($=-\{ e_2, e_1, e_1 \}$). \\
\par
{\bf Example 3.5.} A 4-dimensional (nilpotent) complex left Leibniz algebra $(L, \cdot )$ is given by \\
\par
$e_1 \cdot e_1 = e_4$, $e_1 \cdot e_2 = e_3$, $e_1 \cdot e_3 = e_4$,
\par
$e_2 \cdot e_1 = -e_3 + e_4$,
\par
$e_3 \cdot e_1 = -e_4$\\
\\
(see \cite{Omirov}, Theorem 3.2, algebra $\mathfrak{R}_8$). A linear self-map $\alpha$ of $L$ defined by\\
\par
$\alpha (e_1) = e_1 + e_3 + e_4$,
\par
$\alpha (e_2) = e_2 + e_4$,
\par
$\alpha (e_3) = e_3$,
\par
$\alpha (e_4) = e_4$,\\
\\
is verified to be an endomorphism of $(L, \cdot )$. Again, as in the examples above, define on $L$ an operation ``$\dot{\alpha}$'' by \\
\par
$e_1 \; \dot{\alpha} \; e_1 = e_4$, $e_1 \; \dot{\alpha} \; e_2 = e_3$, $e_1 \; \dot{\alpha} \; e_3 = e_4$,
\par
$e_2 \; \dot{\alpha} \; e_1 = -e_3 + e_4$,
\par
$e_3 \; \dot{\alpha} \; e_1 = -e_4$\\
\\
and then $L_{\alpha} := (L, \dot{\alpha}, \alpha)$ is a multiplicative left Hom-Leibniz algebra. So the theorem implies that the Hom-LY algebra $(L, [,] , \{ , , \}, \alpha)$ constructed from $L_{\alpha}$ is defined by\\
\par
$[e_1 , e_2] = 2e_3 - e_4 \;$ ($=-[e_2, e_1]$), $\;[e_1 , e_3] = 2e_4 \;$ ($=-[e_3, e_1]$),
\par
$\{ e_1, e_2, e_1 \} = e_4 \;$ ($=-\{ e_2, e_1, e_1 \}$).

\end{document}